\documentclass[12pt]{elsarticle}

\usepackage{amsmath} 
\usepackage{amssymb}
\usepackage{amsthm}

\newcommand{\D}{\mathbb{D}}

\newcommand{\I}{\mathbb{I}}
\newcommand{\R}{\mathbb{R}}
\newcommand{\Com}{\mathbb{C}}
\newcommand{\N}{\mathbb{N}}
\newcommand{\nat}{\mathbb{N}}
\newtheorem{Definition}{Definition}
\newtheorem{Theorem}{Theorem}
\newtheorem{Remark}{Remark}
\newtheorem{Example}{Example}

\begin{document}

\begin{frontmatter}



\title{Operational Calculus for the 1st Level General Fractional Derivatives and its Applications}


\author[inst1]{Maryam Alkandari}
\ead{maryam.alkandari@ku.edu.kw}
\affiliation[inst1]{organization={Department of Mathematics, Kuwait University},
            city={Kuwait City},
            postcode={12037}, 
            country={Kuwait}}

\author[inst2]{Yuri Luchko}
\ead{luchko@bht-berlin.de}
\affiliation[inst2]{organization={Department of  Mathematics, Physics, and Chemistry, Berlin University of Applied Sciences and Technology},
            city={Berlin},
            postcode={10587}, 
            country={Germany}}

\begin{abstract}
The 1st level General Fractional Derivatives (GFDs) combine in one definition the GFDs of the Riemann-Liouville type and the regularized GFDs (or the GFDs of the Caputo type) that have been recently introduced and actively studied in the Fractional Calculus literature. In this paper, we first construct  an operational calculus of Mikusi\'nski type for the 1st level GFDs. In particular, it includes the operational calculi for the GFDs of the Riemann-Liouville type and for the regularized GFDs as its particular cases.  In the second part of the paper, this calculus is applied for derivation of the closed form solution formulas to the initial-value problems for the linear fractional differential equations with the 1st level GFDs.
\end{abstract}




\end{frontmatter}


\section{Introduction}
\label{sec1}

Within the last few years, a lot of research in Fractional Calculus (FC) was directed towards the so-called General Fractional Derivatives (GFDs) that are a far reaching generalization of the conventional time-fractional derivatives in form of the integro-differential operators of the Laplace convolution type with the Sonin kernels, see, e.g., \cite{Koch11}-\cite{Tar}.  Along with development of mathematical theory of the GFDs, several important applications of the GFDs have been already suggested, see,  e.g., \cite{ata,baz} for models in linear viscoelasticity involving the GFDs, \cite{mis1,mis2} for applications of the GFDs in anomalous diffusion,  and \cite{Tar1}-\cite{Tar7} for formulation of some general non-local physical models in terms of the GFDs. 

As in the case of the Riemann-Liouville and Caputo fractional derivatives, most of the properties of the GFDs are very different for the GFDs of the Riemann-Liouville type and for the regularized GFDs (or the GFDs of the Caputo type) and thus in some publications only the case of the GFDs of the Riemann-Liouville type was  investigated whereas in other publications only the regularized GFDs were considered. The publications devoted to the fractional differential equations with the GFDs follow exactly  same pattern, i.e., as a rule they present results either only for equations with the regularized GFDs (\cite{Koch11},\cite{Sin18}-\cite{LucYam20}) or only for equations with the GFDs of the Riemann-Liouville type (\cite{Luc21d},\cite{Luc22a}). 

To avoid duplication of research efforts, a concept of the 1st level GFDs was recently introduced in \cite{Luc22c} and \cite{Luc23}. The 1st level GFDs combine in one definition both the GFDs of the Riemann-Liouville type  and the regularized GFD in exactly same way as the generalized Riemann-Liouville fractional derivative (nowadays referred to as the Hilfer fractional derivative) combines definitions of the Riemann-Liouville and the Caputo fractional derivatives, see, e.g., \cite{Hil00}, \cite{Hil19}. This means that any result derived for the  1st level GFDs is automatically valid both for the GFDs of the Riemann-Liouville type and for the regularized GFDs and in particular for the Riemann-Liouville and Caputo fractional derivatives. 

The main focus of this paper is on development of an operational calculus  of Mikusi\'nski type for the 1st level GFDs and on application of this calculus for derivation of the closed form formulas in terms of the so-called convolution series for solutions to  a class of the initial-value problems for the linear fractional differential equations with the the 1st level GFDs.   

The first operational calculus based on the purely algebraic techniques was developed in the 1950s  by Polish mathematician Jan Mikusi\'nski for the first order derivative (\cite{Mik59}, \cite{Yos}). In the framework of  Mikusi\'nski's approach, the first order derivative of a continuous function was interpreted as a multiplication in a special field of convolution quotients (see Section \ref{sec4} for details). Later on, the Mikusi\'nski scheme was extended first to several particular cases of the hyper-Bessel differential operator (\cite{Dit57,Dit63,Mel60}) and then for the general hyper-Bessel differential operator of the $n$-th order (\cite{Dim66}).  In the 1990s, one started to develop  operational calculi of Mikusi\'nski type  for different fractional derivatives including the multiple  Erd\'elyi-Kober fractional derivative (\cite{LucYak94}), the Riemann-Liouville fractional derivative (\cite{HadLuc,LucSri95}, and the  Caputo fractional derivative (\cite{LucGor99}). An operational calculus for the   Hilfer fractional derivative was worked out in \cite{HLT09}.  In \cite{Arr21a,Arr21b},  the operational calculi of Mikusi\'nski type for the Riemann-Liouville and Caputo fractional derivatives with respect to another function were suggested. A Mikusi\'nski type operational calculus for the Prabhakar fractional derivative was introduced in \cite{Arr22a,Arr22b}.

Very recently, the  Mikusi\'nski type operational calculi were developed both for the GFDs of the Riemann-Liouville type  (\cite{Luc22a}) and  for the regularized GFDs (\cite{Luc21c,KHL}). Moreover, in the papers mentioned above, these calculi were applied for derivation of the closed form solution formulas  for the initial-value problems for the multi-term fractional differential equations with the sequential GFDs and the sequential regularized GFDs, respectively. In this paper, we generalize these results for the case of the 1st level GFDs that include the GFDs of the Riemann-Liouville type  and  the regularized GFDs as its particular cases. 

The rest of this paper is organized as follows: In Section \ref{sec2}, we remind the readers of the definitions and some important properties of the General Fractional Integrals (GFIs) and the GFDs with the Sonin kernels including the GFDs of the Riemann-Liouville type, the regularized GFDs, and the 1st level GFDs. Section GFDs \ref{sec3} is devoted to the n-fold GFIs and the n-fold sequential 1st level GFDs. In particular, for the first time in the FC literature, we formulate and prove the first and the second fundamental theorems of FC for the $n$-fold GFIs and the $n$-fold sequential 1st level GFDs. In {Section }\ref{sec4}, we develop  an operational calculus of the Mikusi\'nski type for the 1st level GFDs. In the framework of this operational calculus, the 1st level GFDs and the $n$-fold sequential 1st level  GFDs  are interpreted as multiplication with certain elements of the corresponding fields of convolution quotients. In Section \ref{sec5}, the  operational calculus for the 1st level GFDs is applied for derivations of the closed form solution formulas for the initial-value problems for the fractional differential equations with the 1st level GFDs and for the multi-term fractional differential equations with the $n$-fold sequential 1st level  GFDs. The solutions are provided in terms of the convolution series that are a far reaching generalization of the power law series. 
\section{The 1st level GFDs}
\label{sec2}

For the first time, the 1st level GFDs were introduced in the recent paper \cite{Luc22c} (see also \cite{Luc23} for a generalization of the 1st level GFDs to the case of arbitrary order). In this section, we present a definition of the 1st level GFDs and some of their properties needed for the further discussions, see \cite{Luc22c}, \cite{Luc23} for more details and proofs. 

We start with the definitions of the Sonin kernels and the GFIs and the GFDs with the Sonin kernels. 

\begin{Definition}
\label{d_Son}
A function $\kappa: \R_+ \to \R$ is called a Sonin kernel if there exists a function $k : \R_+ \to \R$ such that the Sonin condition
\begin{equation}
\label{Son}
(\kappa \, *\, k )(t) = \int_0^t \, \kappa(t-\tau)\, k(\tau)\, d\tau\, = \, 1,\ t>0
\end{equation}
holds true. The function $k$ is then called the Sonin kernel associated to the  kernel $\kappa$.
\end{Definition}

The first and for sure the most used and important  pair of the Sonin kernels was considered about two hundreds years ago by Abel in \cite{Abel1}, \cite{Abel2}:
\begin{equation}
\label{power}
\kappa(t) = h_\alpha(t),\ \ k(t) = h_{1-\alpha}(t),\ 0<\alpha < 1,
\end{equation}
where $h_\alpha$ is a power law function defined by the formula
\begin{equation}
\label{h}
 h_\alpha(t):=\frac{t^{\alpha-1}}{\Gamma(\alpha)},\ \alpha >0.
\end{equation}

In the recent paper \cite{Luc_Son}, a general class of the Sonin kernels was introduced as follows:
\begin{equation}
\label{kappa}
\kappa(t) = t^{\alpha-1}\cdot \kappa_1(\lambda t^\beta),\ \kappa_1(t)=\sum_{n=0}^{+\infty}\, a_n t^n, \ a_0 \not = 0,
\end{equation}
\begin{equation}
\label{k}
k(t) = t^{-\alpha} \cdot k_1(\lambda t^\beta),\ k_1(t)=\sum_{n=0}^{+\infty}\, b_n t^n,
\end{equation}
where 
\begin{equation}
\label{cond_s}
0<\alpha < 1,\ \beta >0,\ \lambda \in \R
\end{equation}
and  the coefficients $a_n,\, b_n,\ n\in \nat_0$ of the analytic functions $\kappa_1$ and $k_1$ satisfy the triangular system of the linear equations
\begin{equation}
\label{sys}
\Gamma(\alpha)\Gamma(1-\alpha)a_0 b_{0}=1,\ \sum_{n=0}^N\Gamma(\beta n + \alpha) \Gamma(\beta(N-n)+1-\alpha) a_n b_{N-n}= 0,\ N\in \nat .
\end{equation}

It is worth mentioning that the kernels in form \eqref{kappa} and \eqref{k} with the parameter values $\beta = 1$ and $\lambda =1$ were known already to Sonin. In particular, in \cite{Son},   he derived the famous pair of the Sonin kernels of this kind 
\begin{equation}
\label{Bess}
\kappa(t) = (\sqrt{t})^{\alpha-1}J_{\alpha-1}(2\sqrt{t}),\
k(t) = (\sqrt{t})^{-\alpha}I_{-\alpha}(2\sqrt{t}),\ 0<\alpha <1,
\end{equation}
where $J_\nu$ and $I_\nu$ are the Bessel and the modified Bessel functions, respectively.

An important example of the Sonin kernels in form \eqref{kappa}, \eqref{k} is provided in terms of the three-parameters Mittag-Leffler or the Prabhakar function  (\cite{Prab_2,Giu}):
\begin{equation}
\label{kappa_2}
\kappa(t) = t^{\alpha-1}\, \sum_{n=0}^{+\infty} \frac{(-1)^n (\gamma)_n}{n!\Gamma(\beta n+\alpha) } (\lambda t^\beta)^n  = 
t^{\alpha-1}\, E_{\beta,\alpha}^\gamma(-\lambda t^\beta),
\end{equation}
\begin{equation}
\label{k_2}
k(t) = t^{-\alpha}\, \sum_{n=0}^{+\infty} \frac{(-1)^n (-\gamma)_n}{n!\Gamma(\beta n + 1-\alpha) } (\lambda t^\beta)^n  = 
t^{-\alpha}\, E_{\beta,1-\alpha}^{-\gamma}(-\lambda t^\beta),
\end{equation}
where $\alpha\in (0,\, 1),\ \beta >0,\ \lambda \in \R$ and the function $E_{\beta,\alpha}^\gamma(z)$ is defined by the convergent series                    	
\begin{equation}
\label{Prab}
E_{\beta,\alpha}^\gamma(z) := \sum_{n=0}^{+\infty} \frac{(\gamma)_n}{n!\Gamma(\beta n+\alpha) }\, z^n,\ z,\alpha,\gamma  \in \Com,\ \beta>0.
\end{equation} 

For other examples of the Sonin kernels in terms of the elementary and special functions see, e.g., \cite{Koch11,Han20,Luc21a,Luc21c,Son,Luc_Son,Sam}.

In \cite{Koch11}, for the first time, the integral and integro-differential operators of the Laplace convolution type with the Sonin kernels  $\kappa,\, k$ from a special class $\mathcal{K}$ of kernels were interpreted as the GFI $\I_{(\kappa)}$ and the GFD $\D_{(k)}$ (of the Riemann-Liouville type) and the regularized GFD ${ }_*\D_{(k)}$ (of the Caputo type), respectively:
\begin{equation}
\label{GFI}
(\I_{(\kappa)}\, f)(t) := (\kappa\, *\, f)(t) \, = \, \int_0^t \kappa(t-\tau)f(\tau)\, d\tau,\ t>0,
\end{equation}
\begin{equation}
\label{FDR-L}
(\D_{(k)}\, f)(t) := \frac{d}{dt} (k\, *\, f)(t) = \frac{d}{dt}\, (\I_{(k)}\, f)(t),\ t>0,
\end{equation}
\begin{equation}
\label{FDC}
( _*\D_{(k)}\, f) (t) :=  (\D_{(k)}\, f) (t) - f(0)k(t),\ t>0.
\end{equation}

For a function $f$ that satisfies $f^\prime\in L^1_{loc}(\R_+)$, the regularized GFD  can be represented in the form
\begin{equation}
\label{FDC_1}
( _*\D_{(k)}\, f) (t) =  (\I_{(k)}\, f^\prime)(t) ,\ t>0.
\end{equation}

In \cite{Luc21a,Luc21b,Luc22con,Luc21c}, a theory of the GFIs and the GFDs with the Sonin kernels that belong to the space of functions
\begin{equation}
\label{C-1}
C_{-1}(0,+\infty)\, := \, \{f:\ f(t)=t^{p}f_1(t),\ t>0,\ p > -1,\ f_1\in C[0,+\infty)\}
\end{equation}
was constructed. The class of such kernels is denoted by $\mathcal{L}_1$. 

The basic properties of the GFIs \eqref{GFI} with the kernels from the class $\mathcal{L}_1$ on the space $C_{-1}(0,+\infty)$ are as follows (\cite{Luc21a}):
\begin{equation}
\label{GFI-map}
\I_{(\kappa)}:\, C_{-1}(0,+\infty)\, \rightarrow C_{-1}(0,+\infty),
\end{equation}
\begin{equation}
\label{GFI-com}
\I_{(\kappa_1)}\, \I_{(\kappa_2)} = \I_{(\kappa_2)}\, \I_{(\kappa_1)},
\end{equation}
\begin{equation}
\label{GFI-index}
\I_{(\kappa_1)}\, \I_{(\kappa_2)} = \I_{(\kappa_1*\kappa_2)}.
\end{equation}

As mentioned already by Sonin, the condition \eqref{Son} posed on the kernels $\kappa$ and $k$ of the GFI \eqref{GFI} and the GFD \eqref{FDR-L} or the regularized GFD \eqref{FDC}  ensures that these GFDs are the left-inverse operators to the corresponding GFI on some suitable nontrivial spaces of functions (see, e.g. \cite{Luc21a,Luc21c} for details and proofs). 

For the Sonin kernels $\kappa(t) =h_\alpha (t),\ k(t) = h_{1-\alpha}(t),\ 0< \alpha <1$, the GFI \eqref{GFI} and the GFDs \eqref{FDR-L} and \eqref{FDC} (or \eqref{FDC_1}) define the Riemann-Liouville fractional integral and the Riemann-Liouville and Caputo fractional derivatives of the order $\alpha,\ 0< \alpha <1$, respectively:
\begin{equation}
\label{RLI}
(I^\alpha_{0+}\, f)(t) := (h_\alpha\, * \, f)(t) \, = \, \frac{1}{\Gamma(\alpha)}\, \int_0^t (t-\tau)^{\alpha-1}f(\tau)\, d\tau,\ t>0,
\end{equation}
\begin{equation}
\label{RLD}
(D^\alpha_{0+}\, f)(t) := \frac{d}{dt}\, (I^{1-\alpha}_{0+}\, f)(t),\ t>0,
\end{equation}
\begin{equation}
\label{CD}
( _*D^\alpha_{0+}\, f)(t) :=  (I^{1-\alpha}_{0+}\, f^\prime)(t) \, = \, (D^\alpha_{0+}\, f)(t) - f(0)h_{1-\alpha}(t),\ t>0.
\end{equation}

Because $(I^\alpha_{0+}\, f)(t) \to f(t)$ as $\alpha \to 0$, say, in the sense of the $L_p$-norm for any $p\ge 1$ (see Theorem 2.6 in \cite{SKM}), a natural and standard extension of the definition \eqref{RLI} to the case $\alpha = 0$ is as follows:
\begin{equation}
\label{RLI_0}
(I^0_{0+}\, f)(t) := f(t),\  t>0.
\end{equation}



From the viewpoint of the GFIs and the GFDs with the Sonin kernels, the definition \eqref{RLI_0} can be interpreted as an extension of the Sonin condition \eqref{Son} for the power law kernels \eqref{power} with $0<\alpha < 1$ to the case $\alpha = 0$, i.e., as the formula
\begin{equation}
\label{h_con_0}
(h_0\, *\ h_1)(t)  = (h_1\, *\ h_0)(t)  =  h_1(t) = 1,\ t>0
\end{equation}
that  has to be understood in the sense of the generalized functions (the function $h_0$ plays the role of the Dirac $\delta$-function). 

In its turn, the relation \eqref{h_con_0} allows an extension of the definition of the GFIs to the case of the kernel $h_0$  that is interpreted as the Sonin kernel in the generalized sense:
\begin{equation}
\label{GFI_0}
(\I_{(h_0)}\, f)(t) := (I^0_{0+}\, f)(t) = f(t),\  t>0.
\end{equation}

The 1st level GFDs that we deal with in this paper are a far reaching generalization of the Hilfer fractional derivative of the order $\alpha,\ 0<\alpha <1$ and type $\beta,\ 0\le \beta \le 1$ defined by the formula
\begin{equation}
\label{Hil}
(D^{\alpha,\beta}_{0+}\, f)(t)  = (I_{0+}^{\beta(1-\alpha)}\, \frac{d}{dt}\, I_{0+}^{(1-\alpha)(1-\beta)}\, f)(t),\ \ 
0 < \alpha < 1, \ 0 \le \beta \le 1.
\end{equation}

Properties and applications of the operator \eqref{Hil} are discussed, e.g., in \cite{Hil00,Hil19,HLT09}. In particular, it is worth mentioning that the Hilfer fractional derivative is a particular case of the so-called Djrbashian-Nersessian fractional derivative (\cite{DN,Luc20}). 

Employing the definition \eqref{RLI_0} of the Riemann-Liouville fractional integral of the order $\alpha = 0$,  the Riemann-Liouville and the Caputo fractional derivatives of the order $\alpha \in (0,1)$ can be represented as particular cases of the Hilfer derivative with the type $\beta =0$ and $\beta =1$, respectively:
\begin{equation}
\label{Hil_RL}
(D^{\alpha,0}_{0+}\, f)(x)  = (I_{0+}^{0}\, \frac{d}{dt}\, I_{0+}^{1-\alpha}\, f)(t)\, = \, \frac{d}{dt}\, (I_{0+}^{1-\alpha}\, f)(t)\, = \,(D^\alpha_{0+}\, f)(t),
\end{equation}
\begin{equation}
\label{Hil_C}
(D^{\alpha,1}_{0+}\, f)(t) \, = \,  (I_{0+}^{1-\alpha}\, \frac{d}{dt}\, I_{0+}^{0}\, f)(t)\, = \, (I_{0+}^{1-\alpha}\, \frac{d}{dt}\, f)(t)\, = \,( _*D^\alpha_{0+}\, f)(t).
\end{equation}


Now we introduce the 1st level GFDs that combine in one definition both the GFDs of the Riemann-Liouville type and the regularized GFDs exactly in the same way as the Hilfer derivative combines the Riemann-Liouville and the Caputo fractional derivatives. To this end, the notion of a pair of the Sonin kernels from Definition \ref{d_Son} is first extended to the case of three kernels. 

\begin{Definition}[\cite{Luc22c}]
\label{1l_kernel}
The functions $\kappa,\, k_1,\, k_2:\R_+ \to \R$ that satisfy the condition 
\begin{equation}
\label{1l_cond}
(\kappa\,*\, k_1\, *\, k_2)(t)\, =\, h_1(t)\, = \, 1,\, t>0
\end{equation}
are called the 1st level kernels of the GFDs.
\end{Definition} 

\begin{Remark}
\label{r_kernels}
Because the Laplace convolution is commutative, the order of the functions $\kappa,\, k_1,\, k_2$ in the condition \eqref{1l_cond} can be arbitrary. However, in what follows, we always associate the first function $\kappa$ with the GFI defined by \eqref{GFI} and the functions $k_1$ and $k_2$ with the 1st level  GFD defined by \eqref{1l_GFD}. The function $\kappa$ will be referred to as the kernel associated to the pair of the kernels 
$(k_1,k_2)$. 

It is also worth mentioning that the 1st level kernel is unique as soon as two other kernels from the triple $(\kappa,\, k_1,\, k_2) \in \mathcal{L}_{1}^{1}$ are fixed. This immediately follows from the known fact that the ring $\mathcal{R}_{-1} = (C_{-1}(0,+\infty),+,*)$    does not have any divisors of zero (\cite{Luc21a,Luc21c}). 
\end{Remark}

\begin{Remark}
\label{rSon}
A comparison of the Sonin condition and the condition \eqref{1l_cond} immediately implicates 
that the 1st level kernels $\kappa,\ k_1,\ k_2$ are also the Sonin kernels with the associated kernels $k_1\, *\, k_2$, $\kappa\, * \, k_2$, and $\kappa\, * \, k_1$,  respectively.  
\end{Remark}

In what follows, we denote the set of the 1st level kernels from the space $C_{-1}(0,+\infty)$ 
 by $\mathcal{L}_{1}^{1}$ (1st level kernels of the GFDs of the order less than one) 
and focus on the GFIs and the GFDs with the kernels from the set $\mathcal{L}_{1}^{1}$.

\begin{Example}
Let us consider the power law functions
\begin{equation}
\label{power_ex}
\kappa(t) = h_\alpha(t),\  k_1(t) = h_\gamma(t),\  k_2(t)=h_{1-\alpha-\gamma}(t),\ t>0
\end{equation}
whose parameters satisfy the  conditions
\begin{equation}
\label{power_cond}
0<\alpha<1,\ 0<\gamma < 1-\alpha.
\end{equation}
Under the conditions \eqref{power_cond}, the functions defined by \eqref{power_ex} belong to the space $C_{-1}(0,+\infty)$. They also satisfy the condition \eqref{1l_cond} that immediately follows from the well-known relation
\begin{equation}
\label{alpha_beta}
(h_\alpha\, *\, h_\beta)(t) = h_{\alpha+\beta}(t),\ t>0,\ \alpha>0,\ \beta >0.
\end{equation}
Thus, the functions \eqref{power_ex} are the 1st level kernels  from the class $\mathcal{L}_{1}^{1}$. 
\end{Example}

For a discussion of the properties and other examples of the 1st level kernels as well as the methods for their construction we refer to \cite{Luc22c,Luc23}.

Now we proceed with a definition of the 1st level GFDs with the kernels from the class 
$\mathcal{L}_{1}^{1}$.

\begin{Definition}[\cite{Luc22c}]
\label{1l_GFD_d}
Let $(\kappa,\, k_1,\, k_2) \in \mathcal{L}_{1}^{1}$. 

The 1st level GFD with the kernels $(k_1,\, k_2)$ is defined by the formula
\begin{equation}
\label{1l_GFD}
( _{1L}\D_{(k_1,k_2)}\, f)(t)\, :=\, \left(\I_{(k_1)}\, \frac{d}{dt}\, \I_{(k_2)}\, f\right)(t),
\end{equation}
where $\I_{(k_1)}$ and $\I_{(k_2)}$ are the GFIs defined as in  \eqref{GFI}.
\end{Definition} 


The GFD of the Riemann-Liouville type and the regularized GFD are important particular cases of the 1st level GFD \eqref{1l_GFD} (see the relation \eqref{GFI_0}): 
\begin{equation}
\label{1l_RLD}
( _{1L}\D_{(h_0,k_2)}\, f)(t) = \left(\I_{(h_0)}\, \frac{d}{dt}\, \I_{(k_2)}\, f\right)(t) = \left(\frac{d}{dt}\, \I_{(k_2)}\, f\right)(t) = (\D_{(k_2)}\, f)(t),
\end{equation}
\begin{equation}
\label{1l_CD}
( _{1L}\D_{(k_1,h_0)}\, f)(t) = \left(\I_{(k_1)}\, \frac{d}{dt}\, \I_{(h_0)}\, f\right)(t) = \left(\I_{(k_1)} f^\prime \right)(t) = ( _*\D_{(k_1)}\, f)(t).
\end{equation}

Another significant particular case of the 1st level GFD \eqref{1l_GFD} is the Hilfer fractional derivative with the power law kernels  as in \eqref{power_ex}:
\begin{equation}
\label{1l_GFD_Hilfer}
( _{1L}\D_{(h_\gamma,h_{1-\alpha-\gamma})}\, f)(t)\, =\, \left(I_{0+}^{\gamma}\, \frac{d}{dt}\, I_{0+}^{1-\alpha-\gamma}\, f\right)(t).
\end{equation}
Please note that the operator \eqref{1l_GFD_Hilfer} defines the Hilfer fractional derivative  in a slightly different parametrization compared to the one used in the formula \eqref{Hil}.

For the 1st level kernels $(\kappa,\, k_1,\, k_2) \in \mathcal{L}_{1}^{1}$, the GFI with the kernel $\kappa$ and the 1st level GFD generated by the kernels $(k_1,\, k_2)$ build a calculus in the sense that the GFI and the GFDs satisfy two fundamental theorems of calculus, of course, in a generalized form. For the first time, these results were derived in \cite{Luc22c} for the 1st level GFDs that we deal with in this paper and  in \cite{Luc23} for the 1st level GFDs of arbitrary order. 

\begin{Theorem}[\cite{Luc22c}]
\label{t-FT1}
Let $(\kappa,\, k_1,\, k_2) \in \mathcal{L}_{1}^{1}$.

Then the 1st level GFD \eqref{1l_GFD} is a left-inverse operator to the GFI \eqref{GFI}  on the space of functions $\I_{(k_1)}(C_{-1}(0,+\infty))$,  i.e., the relation 
\begin{equation}
\label{FT1L}
( _{1L}\D_{(k_1,k_2)}\, \I_{(\kappa)}\, f) (t) = f(t),\ t>0,\  f\in \I_{(k_1)}(C_{-1}(0,+\infty))
\end{equation}
holds valid, where
\begin{equation}
\label{C_k}
\I_{(k)}(C_{-1}(0,+\infty)) := \{f:\ f(t)=(\I_{(k)}\, \phi)(t),\ \phi\in C_{-1}(0,+\infty)\}.
\end{equation}
\end{Theorem}

As in the case of the conventional calculus,  the 2nd fundamental theorem of FC for the 1st level GFD is formulated  on an essentially narrower space of functions compared to the one employed in Theorem \ref{t-FT1}, namely on the space 
\begin{equation}
\label{C1(k)}
C_{-1,(k)}^{1}(0,+\infty)= \{ f\in C_{-1}(0,+\infty):\, (\D_{(k)}\, f)(t) \in C_{-1}(0,+\infty) \}.
\end{equation}

\begin{Theorem}[\cite{Luc22c}]
\label{t-FT2}
Let $(\kappa,\, k_1,\, k_2) \in \mathcal{L}_{1}^{1}$.

Then the formula
\begin{equation}
\label{FT2L}
(\I_{(\kappa)}\,  _{1L}\D_{(k_1,k_2)}\, f) (t) = f(t)\, - \, (\I_{(k_2)}\, f)(0)\, (\kappa\, *\, k_1)(t),\  t>0
\end{equation}
holds valid for any function $f\in C_{-1,(k_2)}^{1}(0,+\infty).$
\end{Theorem}

\begin{Remark}
\label{r2_1}
According to the formula \eqref{1l_GFD_Hilfer}, the Hilfer fractional derivative is a particular case of the 1st level GFD. For this derivative, the formula \eqref{FT2L} takes the well-known form (\cite{Hil00,HLT09}):
\begin{equation}
\label{proj_H}
(I_{0+}^{\alpha}\, _{1L}\D_{(h_\gamma,h_{1-\alpha-\gamma})}\, f)(t) = f(t) - (I_{0+}^{1-\alpha-\gamma}\, f)(0)\, h_{\alpha+\gamma}(t),\ t>0.
\end{equation}
\end{Remark}

\begin{Remark}
\label{r2}
It is worth mentioning that the formula \eqref{FT2L} 
can be rewritten in terms of the projector operator of the 1st level GFD:
\begin{equation}
\label{proj_1}
(P_{1L}\, f)(t) := f(t) - (\I_{(\kappa)}\,  _{1L}\D_{(k_1,k_2)}\, f) (t) =   (\I_{(k_2)}\, f)(0)\, (\kappa\, *\, k_1)(t),\ t>0.
\end{equation}
In its turn,  the formula \eqref{proj_1} for the projector operator determines the natural initial conditions while dealing with the fractional differential equations that contain the 1st level GFDs. In Section \ref{sec5}, we consider some equations with the initial conditions of this kind. 
\end{Remark}


\section{The n-fold GFIs and the n-fold sequential 1st level GFDs}
\label{sec3}

In  this section, we define the $n$-fold GFIs and the $n$-fold sequential 1st level GFDs  that are new objects not yet considered in the literature and discuss some of their properties. 

\begin{Definition}
\label{d_n}
Let $(\kappa,\, k_1,\, k_2) \in \mathcal{L}_{1}^{1}$ and $n \in \N$. 

The $n$-fold  GFIs and the $n$-fold  sequential 1st level GFDs are defined as follows, respectively:
\begin{equation}
\label{GFIn}
(\I_{(\kappa)}^{<n>}\, f)(t) := (\underbrace{\I_{(\kappa)} \ldots \I_{(\kappa)}}_{n\ \mbox{times}}\, f)(t),\  t>0, 
\end{equation}
\begin{equation}
\label{GFDLn}
( _{1L}\D_{(k_1,k_2)}^{<n>}\, f)(t)
:= (\underbrace{ _{1L}\D_{(k_1,k_2)} \ldots  _{1L}\D_{(k_1,k_2)}}_{n\ \mbox{times}}\, f)(t),\  t>0.
\end{equation}
\end{Definition}

\begin{Remark}
\label{rn=0}
For $n=0$, the $n$-fold  GFI and the $n$-fold  sequential 1st level GFD are interpreted as an identity operator:
\begin{equation}
\label{GFIn=0}
(\I_{(\kappa)}^{<0>}\, f)(t) := f(t),\ ( _{1L}\D_{(k_1,k_2)}^{<0>}\, f)(t):= f(t),\ t>0. 
\end{equation}
\end{Remark}

Using the well-known properties of the Laplace convolution, the $n$-fold  GFI can be represented as an integral operator of the Laplace convolution type:
\begin{equation}
\label{GFIn-1}
(\I_{(\kappa)}^{<n>}\, f)(t) = (\kappa^{<n>}\, *\, f)(t),\  t>0,
\end{equation}
where $\kappa^{<n>}$ stands for a convolution power
\begin{equation}
\label{conpower}
\kappa^{<n>}(t):= \begin{cases}
\kappa(t), & n=1,\\
(\underbrace{\kappa* \ldots * \kappa}_{n\ \mbox{times}})(t),& n=2,3,\dots .
\end{cases}
\end{equation}

\begin{Remark}
\label{rcpn=0}
The formula \eqref{GFIn-1} is valid also for the $n$-fold GFI with $n=0$ defined by \eqref{GFIn=0} provided that the definition\eqref{conpower} of the convolution power is extended to the case $n=0$ as follows:
\begin{equation}
\label{cpn=0}
\kappa^{<0>}(t):=\delta(t),
\end{equation}
where $\delta$ is the Dirac $\delta$-function. 
\end{Remark}

Please note that even if the functions $\kappa^{<n>}(t),\ n\in \N$ belong to the space $C_{-1}(0,+\infty)$,  for $n\ge 2$, they are not always the Sonin kernels  and thus the operator defined by  \eqref{GFIn} or \eqref{GFIn-1} is not always a GFI, see \cite{Luc21c} for details.  

For the first time, the $n$-fold sequential GFDs of the Riemann-Liouville type and the $n$-fold sequential regularized GFDs were introduced and investigated  in \cite{Luc22con}. They are particular cases of the $n$-fold sequential  1st level GFD defined by \eqref{GFDLn}, see the formulas \eqref{1l_RLD} and \eqref{1l_CD}. 

In its turn, the $n$-fold sequential GFDs of the Riemann-Liouville type and the $n$-fold sequential regularized GFDs are generalizations of the well-known Riemann-Liouville and Caputo sequential fractional derivatives.

In the rest of this section, we derive some important properties of the $n$-fold  GFI and the $n$-fold  sequential 1st level GFD. 

Repeatedly applying Theorem \ref{t-FT1}, we arrive at the following formulation of the first fundamental theorem of FC for the $n$-fold sequential  1st level  GFDs:

\begin{Theorem}
\label{t3-n}
Let $(\kappa,\, k_1,\, k_2) \in \mathcal{L}_{1}^{1}$ and $n\in \nat$.

Then,  the $n$-fold sequential 1st level GFD \eqref{GFDLn}   is a left-inverse operator to the $n$-fold  GFI \eqref{GFIn} on the space $\I_{(k_1)}(C_{-1}(0,+\infty))$: 
\begin{equation}
\label{FT1L-n}
( _{1L}\D_{(k_1,k_2)}^{<n>}\, \I_{(\kappa)}^{<n>}\, f) (t) = f(t),\ t>0,\  f\in \I_{(k_1)}(C_{-1}(0,+\infty)).
\end{equation}
\end{Theorem}

Indeed, Theorem \ref{t-FT1} implicates the recurrent formula ($n\in \nat$ and $n\ge 2$)
$$
( _{1L}\D_{(k_1,k_2)}^{<n>}\, \I_{(\kappa)}^{<n>}\, f) (t) = ( _{1L}\D_{(k_1,k_2)}^{<n-1>}\,(  _{1L}\D_{(k_1,k_2)}\, \I_{(\kappa)}\, ( \I_{(\kappa)}^{<n-1>}\, f))) (t) =
$$
$$
( _{1L}\D_{(k_1,k_2)}^{<n-1>}\, \I_{(\kappa)}^{<n-1>}\, f) (t) 
$$
and the statement of Theorem \ref{t3-n} follows immediately. 

As to the second fundamental theorem of FC for the $n$-fold GFIs and the $n$-fold sequential 1st level GFDs, we formulate and prove it for the functions from the space 

\begin{eqnarray}
 \label{space}
C_{-1,(k_1,k_2)}^{m}(0,+\infty) && =  \{ f:\   _{1L}\D_{(k_1,k_2)}^{<j>}\, f \in C_{-1}(0,+\infty)\ \wedge \\
   &&  \D_{(k_2)}\,  _{1L}\D_{(k_1,k_2)}^{<j>}\, f \in C_{-1}(0,+\infty),\ j=0,\dots,m\}. \nonumber
\end{eqnarray}

\begin{Theorem}
\label{tsft-n}
Let $(\kappa,\, k_1,\, k_2) \in \mathcal{L}_{1}^{1}$ and $n\in \nat$.

Then the formula
\begin{equation}
\label{sFTLn}
(\I_{(\kappa)}^{<n>}\,  _{1L}\D_{(k_1,k_2)}^{<n>}\, f) (t) = f(t) - \sum_{j=0}^{n-1}(\I_{(k_2)}\,  _{1L}\D_{(k_1,k_2)}^{<j>}\, f)(0)(\kappa^{<j+1>}\, *\, k_1) (t)
\end{equation}
holds valid for the functions $f$ from the space $C_{-1,(k_1,k_2)}^{n-1}(0,+\infty)$ defined as in \eqref{space}. 
\end{Theorem}

\begin{proof}
The main ingredient of the proof is  the recurrent formula ($n\in \nat$ and $n\ge 2$)
$$
( \I_{(\kappa)}^{<n>}\, _{1L}\D_{(k_1,k_2)}^{<n>}\,  f) (t) = (\I_{(\kappa)}^{<n-1>}\,  (\I_{(\kappa)} \,  _{1L}\D_{(k_1,k_2)}\,( _{1L}\D_{(k_1,k_2)}^{<n-1>}\, f)))(t) =
$$ 
$$
(\I_{(\kappa)}^{<n-1>}\, (   _{1L}\D_{(k_1,k_2)}^{<n-1>}\, f)(\cdot) - (\I_{(k_2)}\,  _{1L}\D_{(k_1,k_2)}^{<n-1>}\, f)(0)\, (\kappa\, *\, k_1)(\cdot))(t)= 
$$
$$
(\I_{(\kappa)}^{<n-1>}\,    _{1L}\D_{(k_1,k_2)}^{<n-1>}\, f)(t) - (\I_{(k_2)}\,  _{1L}\D_{(k_1,k_2)}^{<n-1>}\, f)(0)\, (\kappa^{<n>}\, *\, k_1)(t) 
$$
that follows from Theorem \ref{t-FT2} under the condition that both $ _{1L} \, \D_{(k_1,k_2)}^{<n-1>}\, f$ and $\D_{(k_2)} \, _{1L} \, \D_{(k_1,k_2)}^{<n-1>}\, f$ are from the space $C_{-1}(0,+\infty)$.

Because of the inclusion  $f\in C_{-1,(k_1,k_2)}^{n-1}(0,+\infty)$, the functions $ _{1L} \, \D_{(k_1,k_2)}^{<j>}\, f $ and $\D_{(k_2)} \, _{1L} \, \D_{(k_1,k_2)}^{<j>}\, f$ belong to the space $C_{-1}(0,+\infty)$ for all $j=0,1,\dots n-1$ and thus the recurrent formula can be applied $n$ times that leads to the formula \eqref{sFTLn} and the proof is completed.
\end{proof}

\begin{Remark}
\label{rn=1}
For $n=1$, the result formulated in Theorem \ref{tsft-n} coincides with the second fundamental theorem for the 1st level GFD (Theorem \ref{t-FT2}). In this case, we interpret the space $C_{-1,(k_1,k_2)}^{0}(0,+\infty)$ as the space $ C_{-1,(k_2)}^{1}(0,+\infty)$ defined by \eqref{C1(k)} (see the formula \eqref{GFIn=0} for definition of the $n$-fold  sequential 1st level GFD in the case $n=0$). 
\end{Remark}

\begin{Remark}
\label{rprojn}
The formula \eqref{sFTLn} can be rewritten in terms of the  projector operator of the $n$-fold sequential 1st level GFD:
\begin{equation}
\label{proj_n}
(P_{1L}^{<n>}\, f)(t) := f(t) -(\I_{(\kappa)}^{<n>}\,  _{1L}\D_{(k_1,k_2)}^{<n>}\, f) (t) = 
\end{equation}
$$
\sum_{j=0}^{n-1}(\I_{(k_2)}\,  _{1L}\D_{(k_1,k_2)}^{<j>}\, f)(0)(\kappa^{<j+1>}\, *\, k_1) (t).
$$

The coefficients $(\I_{(k_2)}\,  _{1L}\D_{(k_1,k_2)}^{<j>}\, f)(0),\ j= 0,1,\dots,n-1$ by the functions $(\kappa^{<j+1>}\, *\, k_1) (t)$ in the sum at the right-hand side of the formula \eqref{proj_n} determine the natural initial conditions while dealing with the fractional differential equations that contain the $n$-fold sequential 1st level GFD. 
\end{Remark}

The two important examples of the formula \eqref{sFTLn} are its particular cases for the $n$-fold sequential GFDs of the Riemann-Liouville type and the $n$-fold sequential regularized GFDs. 

The $n$-fold sequential GFD of the Riemann-Liouville type is a particular case of the $n$-fold sequential 1st level GFD with the kernel $k_1=h_0$. Let us denote the kernel $k_2$ by  $k$. Then the kernels $\kappa$ and $k$ are the Sonin kernels and we arrive at the following result (\cite{Luc22con}): 
\begin{equation}
\label{sFTLn_RL}
(\I_{(\kappa)}^{<n>}\, \D_{(k)}^{<n>}\, f) (t) = 
f(t) - \sum_{j=0}^{n-1}\left( \I_{(k)}\, \D_{(k)}^{<j>}\, f\right)(0)\kappa^{<j+1>}(t),
\end{equation}
 where $\I_{(\kappa)}^{<n>}$ is the $n$-fold GFI \eqref{GFIn} and $\D_{(k)}^{<n>}$ is the $n$-fold sequential GFD in the Riemann-Liouville sense. 

According to Theorem \ref{tsft-n}, the formula \eqref{sFTLn_RL} is valid for the functions $f$ from the space 
 \begin{equation}
\label{C1(k)_n}
C_{-1,(k)}^{n}(0,+\infty)= \{ f:\, (\D_{(k)}^{<j>}\, f)(t) \in C_{-1}(0,+\infty),\ j=0,1,\dots,n \}.
\end{equation}

Now we proceed with  the second fundamental theorem for the $n$-fold sequential regularized GFDs. It is a particular case of the general result derived for the $n$-fold sequential 1st level GFDs with $k_2=h_0$. Let us denote the kernel $k_1$ by $k$. Then the formula \eqref{sFTLn} takes the following form (\cite{Luc22con}):
\begin{equation}
\label{sFTLn-C}
(\I_{(\kappa)}^{<n>}\,_*\D_{(k)}^{<n>}\, f) (t) = f(t) - f(0) - \sum_{j=1}^{n-1}\left(\,_*\D_{(k)}^{<j>}\, f\right)(0)\left( \kappa^{<j+1>} \, *\,  k \right)(t),
\end{equation}
where $\I_{(\kappa)}^{<n>}$ is the $n$-fold GFI \eqref{GFIn} and $\,_*\D_{(k)}^{<n>}$ is the $n$-fold sequential regularized GFD. 

The formula \eqref{sFTLn-C} holds true for the functions $f$ from the space
 \begin{equation}
\label{C1(k)_n_C}
\,_*C_{-1,(k)}^{n-1}(0,+\infty)= \{ f:\, (\,_*\D_{(k)}^{<j>}\, f)(t) \in C_{-1}(0,+\infty),\ j=0,1,\dots,n-1 \}.
\end{equation}

Finally, we formulate the second fundamental theorem for the $n$-fold sequential Hilfer fractional derivative generated by the right-hand side of the formula \eqref{1l_GFD_Hilfer}. In the case of the Hilfer fractional derivative, the kernels of the corresponding 1st level GFD are the power law functions ($\kappa = h_\alpha$, $k_1 = h_\gamma$, and $k_2 = h_{1-\alpha - \gamma}$) and the formula \eqref{sFTLn} takes the form
\begin{equation}
\label{sFTLn-H}
(I_{0+}^{\alpha\, n}\, _{1L}\D_{(h_\gamma,h_{1-\alpha-\gamma})}^{<n>}\, f) (t) = f(t) - \sum_{j=0}^{n-1}\left(I_{0+}^{1-\alpha - \gamma}\, _{1L}\D_{(h_\gamma,h_{1-\alpha-\gamma})}^{<j>}\, f\right)(0)\, h_{\alpha j +\alpha +\gamma}(t).
\end{equation}


\section{Operational calculus for 1st level GFDs}
\label{sec4}

In this section, we develop a Mikusi\'nski type operational calculus for the 1st level GFD defined as in \eqref{1l_GFD}. 

The main idea behind any operational calculus of Mikusi\'nski type  is  an interpretation of some differential, integral, or integro-differential operators as purely algebraic operations. In the framework of this interpretation, the differential, integral, or integro-differential equations with these operators can be reduced to some algebraic equations. The solutions to these algebraic equations can be interpreted as the generalized solutions to the corresponding differential, integral, or integro-differential equations. By means of the so-called operational relations,  these generalized solutions can be interpreted as solutions in the strong sense. 
The main idea behind any operational calculus of Mikusi\'nski type  is  an interpretation of some differential, integral, or integro-differential operators as purely algebraic operations. In the framework of this interpretation, the differential, integral, or integro-differential equations with these operators can be reduced to some algebraic equations. The solutions to these algebraic equations can be interpreted as the generalized solutions to the corresponding differential, integral, or integro-differential equations. In some cases,  by means of the so-called operational relations, the generalized solutions can be interpreted as solutions in the strong sense.   

The starting point of the Mikusi\'nski type operational calculi for different kinds of the fractional derivatives developed so far is the  statement that  $\mathcal{R}_{-1} = (C_{-1}(0,+\infty),+,*)$     with  the  usual
addition $+$ and  multiplication $*$ in form of  the  Laplace convolution is a commutative ring without divisors of zero  first proved in \cite{LucGor99}. 

By definition,  the GFI \eqref{GFI} can be interpreted as a multiplication with the Sonin kernel $\kappa$ on the ring $\mathcal{R}_{-1}$:
\begin{equation}
\label{GFI_op}
(\I_{(\kappa)}\, f)(t) = (\kappa\, *\, f)(t),\ f\in C_{-1}(0,+\infty),\ t>0. 
\end{equation}

As shown in \cite{LucGor99}, the ring $\mathcal{R}_{-1}$ does not possess a unity element with respect to multiplication. According to Theorem \ref{t-FT1}, the 1st level GFD \eqref{1l_GFD} is a left-inverse operator to the GFI. Thus, a representation of the 1st level GFD as a multiplication on the ring $\mathcal{R}_{-1}$ is not possible. 

Another key element of the Mikusi\'nski type operational calculi for the fractional derivatives is an extension of the ring $\mathcal{R}_{-1}$  to a field of convolution quotients. This construction was first suggested by  Mikusi\'nski in the framework of his operational calculus for the 1st order derivative developed for the ring of functions continuous on the positive semi-axes and then adjusted  to the case of the Riemann-Liouville fractional derivative and the ring  $\mathcal{R}_{-1}$  in \cite{LucSri95}. In what follows, we shortly outline this standard procedure. 

We start with defining an equivalence relation 
$$
(f_1,\, g_1)\sim     (f_2,\, g_2)    \Leftrightarrow     (f_1\, * \,
g_2)(t) = (f_2\, * \,  g_1)(t),\ (f_1,\, g_1),\ (f_2,\, g_2) \in C_{-1}^2
$$
on the set
$$
C_{-1}^2:= C_{-1}(0,+\infty) \times  (C_{-1}(0,+\infty) \setminus \{0\}).
$$

For the sake of convenience, the classes of equivalences $C_{-1}^2/\sim$ will be denoted as quotients:
$$
\frac{f}{g}:= \{ (f_1,\, g_1) \in C_{-1}^2:\ (f_1,\, g_1) \sim (f,\, g) \}.
$$

With these notations, the operations $+$ and $\cdot$ on the set  $C_{-1}^2/\sim$ are defined as in the case of the rational numbers:
$$
\frac{f_1}{g_1} + \frac{f_2}{g_2}:= \frac{ f_1\, *\,  g_2\, +\, f_2\, * g_1}
{g_1\, *\, g_2},
$$
$$
\frac{f_1}{g_1}\cdot  \frac{f_2}{g_2} :=\frac{f_1\, *\,  f_2}{g_1\, *\,  g_2}.
$$

It is worth mentioning that the operations $+$ and $\cdot$ are correctly defined, i.e., their outcomes  do not depend on the representatives of the  equivalence classes that are taken for this operations.

As first shown in \cite{LucGor99}, 
the triple $\mathcal{F}_{-1} = (C_{-1}^2/\sim,\ +,\ \cdot)$ is a  field that is referred to as the  field of convolution quotients.

Evidently, the ring $\mathcal{R}_{-1}$ can be embedded into the field
$\mathcal{F}_{-1}$, say,  by the mapping:
\begin{equation}
\label{emb}
f  \mapsto \frac{f\, *\, \kappa}{\kappa},
\end{equation}
where $\kappa \in C_{-1}(0,+\infty)$ is the Sonin kernel of the GFI \eqref{GFI}.

Because the space $C_{-1}(0,+\infty)$ is a vector space, the set of $C_{-1}^2/\sim$ of  equivalence classes is also a vector space with the scalar multiplication defined as follows:
$$
\lambda\, \frac{f}{g}:=\frac{\lambda\, f}{g},\ \frac{f}{g}\in C_{-1}^2/\sim,\ \lambda \in \R \ 
\mbox{ or } \lambda \in \Com.
$$

Please note that  we have to distinguish between multiplication with the scalar $\lambda$ in the vector
space  $C_{-1}^2/\sim$ and  multiplication with the constant
function $\{ \lambda \}$ in the field $\mathcal{F}_{-1}$:
$$
\{\lambda\}\cdot \frac{f}{g} = \frac{\{\lambda\}\, *\,  f}{g},\ \frac{f}{g}\in \mathcal{F}_{-1}.
$$

In contrast to the ring $\mathcal{R}_{-1}$, the field $\mathcal{F}_{-1}$ possesses  a unity element with respect to multiplication that we denote by $I$:
\begin{equation}
\label{I}
I = \frac{f}{f} = \{ (f,\, f):\, f\in C_{-1}(0,+\infty),\ f\not \equiv 0\}.
\end{equation}
For convenience, in what follows, we associate the unity element $I$ with the representative $(\kappa,\, \kappa)$ of this equivalence class. 

As was shown in \cite{Luc21c}, the unity element $I$ does not belong to the ring $\mathcal{R}_{-1}$ and thus it is a generalized function or the so-called hyper-function, not a conventional one. In our operational calculus, the unity element $I$  plays the role of the Dirac $\delta$-function . 

Now we are going to interpret the 1st level GFD \eqref{1l_GFD} as a multiplication  on the field $\mathcal{F}_{-1}$ of convolution quotients. For this aim,  we first introduce  the following important hyper-function: 

\begin{Definition}[\cite{Luc21c}]
\label{d4}
The element (equivalence class) of the field $\mathcal{F}_{-1}$ with the representative 
\begin{equation}
\label{alg}
S_\kappa := \frac{\kappa}{\kappa^{<2>}}
\end{equation}
 is called algebraic inverse to the GFI \eqref{GFI} with the Sonin kernel $\kappa$.
\end{Definition}


It is easy to see that $S_{\kappa}$ is the  inverse element to  the Sonin kernel $\kappa \in \mathcal{R}_{-1}$ that is interpreted as the element $\frac{\kappa^{<2>}}{\kappa}$ of the field $\mathcal{F}_{-1}$:
\begin{equation}
\label{inverse}
S_\kappa \cdot  {\kappa} = \frac{\kappa}{\kappa^{<2>}} \cdot \frac{\kappa^{<2>}}{\kappa}
= \frac{\kappa^{<3>}}{\kappa^{<3>}} =  I.
\end{equation}

The algebraic inverse to the GFI \eqref{GFI} defined above is now used to introduce a concept of a 1st level algebraic  general fractional derivative (1st level AGFD). 

\begin{Definition}
\label{d5}
Let $(\kappa,\, k_1,\, k_2) \in \mathcal{L}_{1}^{1}$.

The 1st level AGFD is defined by the expression
\begin{equation}
\label{AGFD}
 _{1L}\D_{(k_1,k_2)}\, f := S_\kappa \cdot f - S_\kappa \cdot (P_{1L}\, f)(t),
\end{equation}
where the  projector operator $P_{1L}$ of the 1st level GFD  is given by the formula \eqref{proj_1} and  the function $(P_{1L}\, f)(t)$  at the right-hand side of \eqref{AGFD} is interpreted as an element of the  field $\mathcal{F}_{-1}$.
\end{Definition}

\begin{Remark}
\label{r_AGFD}
The formula \eqref{proj_1} for the projector operator $P_{1L}$ along with the formula \eqref{inverse} implicates another form of the 1st level AGFD:
\begin{equation}
\label{AGFD_1}
 _{1L}\D_{(k_1,k_2)}\, f = S_\kappa \cdot f -  (\I_{(k_2)}\, f)(0)\, k_1(t),
\end{equation}
where the Sonin kernel $k_1$ is interpreted as an element of the field $\mathcal{F}_{-1}$. 
\end{Remark}

As already mentioned, the 1st level GFD \eqref{1l_GFD} does  exist for the functions from the space
$C_{-1,(k_2)}^{1}(0,+\infty)$ defined by \eqref{C1(k)}. For existence of the 1st level AGFD \eqref{AGFD_1}, a much weaker inclusion $\I_{(k_2)}\, f \in C[0,+\infty)$ is sufficient. In this sense, the 1st level AGFD is a generalized derivative that makes sense also for the functions whose conventional 1st level GFD does not exist. However, the 1st level AGFD coincides with the 1st level GFD on the space $C_{-1,(k_2)}^{1}(0,+\infty)$ and thus we employ the same notation for both derivatives.  

\begin{Theorem}
\label{t9}
On the space $C_{-1,(k_2)}^{1}(0,+\infty)$, the 1st level GFD exists in the usual sense and coincides with the 1st level AGFD, i.e., the relation
\begin{equation}
\label{AGFD-f}
( _{1L}\D_{(k_1,k_2)}\, f)(t) ={ }_{1L}\D_{(k_1,k_2)}\, f   = S_\kappa \cdot f -  (\I_{(k_2)}\, f)(0)\, k_1(t)
\end{equation}
hold true provided all functions from the ring $\mathcal{R}_{-1}$ in this relation are interpreted as the elements of the field $\mathcal{F}_{-1}$.
\end{Theorem}

\begin{proof}

The inclusion $f\in C_{-1,(k_2)}^{1}(0,+\infty)$  means that $f\in C_{-1}(0,+\infty)$ and $\D_{(k_2)} f = \frac{d}{dt} \I_{(k_2)} f \in  C_{-1}(0,+\infty)$. According to a result derived in \cite{LucGor99} (see property 3) on page 211), the last inclusion implicates that $\I_{(k_2)}\, f \in C[0,+\infty)$. This proves existence of the 1st level AGFD. Because of the inclusion $\D_{(k_2)} f  \in  C_{-1}(0,+\infty)$, the mapping property \eqref{GFI-map} ensures existence of the 1st level GFD. 

To prove the relation \eqref{AGFD-f}, we embed  the formula \eqref{FT2L} from Theorem \ref{t-FT2} (the second fundamental theorem of FC for the 1st level GFD) into the field $\mathcal{F}_{-1}$ and multiply it by the element $S_\kappa$:
$$
S_\kappa \cdot (\I_{(\kappa)}\,  _{1L}\D_{(k_1,k_2)}\, f) (t) = S_\kappa \cdot (f(t)\, - \, (\I_{(k_2)}\, f)(0)\, (\kappa\, *\, k_1)(t)).
$$
Due to the formulas \eqref{GFI_op} and \eqref{inverse}, the left-hand side of the obtained formula is an embedding of the function $( _{1L}\D_{(k_1,k_2)}\, f)(t)$ into the field $\mathcal{F}_{-1}$ whereas its right-hand side coincides with the definition of the 1st level AGFD that completes the proof of the theorem. 
\end{proof}

Both Definition \ref{d5} and Theorem \ref{t9} can be easily extended to the case of the $n$-fold sequential 1st level GFD defined by \eqref{GFDLn}. 

\begin{Definition}
\label{d6}
Let $(\kappa,\, k_1,\, k_2) \in \mathcal{L}_{1}^{1}$ and $n\in \N$.

The $n$-fold sequential 1st level AGFD is defined by the expression
\begin{equation}
\label{AGFD_n}
 _{1L}\D_{(k_1,k_2)}^{<n>}\, f := S_\kappa^n \cdot f - S_\kappa^n \cdot (P_{1L}^{<n>}\, f)(t),
\end{equation}
where the  projector operator $P_{1L}^{<n>}$ of the $n$-fold sequential 1st level GFD  is given by the formula \eqref{proj_n} and  the function $(P_{1L}^{<n>}\, f)(t)$  at the right-hand side of \eqref{AGFD_n} is interpreted as an element of the  field $\mathcal{F}_{-1}$.
\end{Definition}

\begin{Remark}
\label{r_AGFD_n}
The formula \eqref{proj_n} for the projector operator $P_{1L}^{<n>}$ leads to another form of the $n$-fold sequential 1st level AGFD:
\begin{equation}
\label{AGFD_1_n}
 _{1L}\D_{(k_1,k_2)}^{<n>}\, f = S_\kappa^n \cdot f -  
 \sum_{j=0}^{n-1}(\I_{(k_2)}\,  _{1L}\D_{(k_1,k_2)}^{<j>}\, f)(0)\, S_\kappa^{n-j-1} \cdot k_1(t),
\end{equation}
where the Sonin kernel $k_1$ is interpreted as an element of the field $\mathcal{F}_{-1}$ and $S_\kappa^0$ is defined as the unity element $I$ of $\mathcal{F}_{-1}$ with respect to multiplication operation. 
\end{Remark}

The $n$-fold sequential 1st level AGFD is  a kind of a generalized derivative because its domain is much broader compared to the one of the $n$-fold sequential 1st level GFD.  However, for the functions from the space $C_{-1,(k_1,k_2)}^{n-1}(0,+\infty)$ defined as in \eqref{space},
we have the following result:

\begin{Theorem}
\label{tsft-n-a}
Let $(\kappa,\, k_1,\, k_2) \in \mathcal{L}_{1}^{1}$ and $n\in \nat$.

For a function $f$ from the space $ C_{-1,(k_1,k_2)}^{n-1}(0,+\infty)$ defined as in \eqref{space}, its $n$-fold sequential 1st level GFD \eqref{GFDLn} exists in the usual sense and coincides with the $n$-fold sequential 1st level AGFD \eqref{AGFD_1_n}:
\begin{equation}
\label{AGFD-f_n}
( _{1L}\D_{(k_1,k_2)}^{<n>}\, f)(t) ={ }_{1L}\D_{(k_1,k_2)}^{<n>}\, f   = S_\kappa^n \cdot f - 
 \sum_{j=0}^{n-1}(\I_{(k_2)}\,  _{1L}\D_{(k_1,k_2)}^{<j>}\, f)(0)\, S_\kappa^{n-j-1} \cdot k_1(t)
\end{equation}
provided all functions from the ring $\mathcal{R}_{-1}$ in this formula  are interpreted as the elements of the field $\mathcal{F}_{-1}$.
\end{Theorem}

Indeed, the formula \eqref{AGFD-f_n} immediately follows from the second fundamental theorem of FC for the $n$-fold sequential 1st level GFD. Multiplying the formula \eqref{sFTLn}  embedded into the field $\mathcal{F}_{-1}$ by the element $S_\kappa^n$ we get the relation
$$
S_\kappa^n \, \cdot \, (\I_{(\kappa)}^{<n>}\,  _{1L}\D_{(k_1,k_2)}^{<n>}\, f) (t) = 
$$
$$
S_\kappa^n\, \cdot f  - S_\kappa^n\, \cdot \, (\sum_{j=0}^{n-1}(\I_{(k_2)}\,  _{1L}\D_{(k_1,k_2)}^{<j>}\, f)(0)(\kappa^{<j+1>}\, *\, k_1) (t))
$$
that can be rewritten in form \eqref{AGFD-f_n} by using the formula \eqref{inverse}.

Theorems \ref{t9} and \ref{tsft-n-a}  provide representations of the 1st level GFD and the $n$-fold sequential 1st level GFD in terms of the algebraic operations on the field $\mathcal{F}_{-1}$ of convolution quotients. In particular, one can employ these representations for rewriting the multi-term fractional differential equations with the $n$-fold sequential 1st level GFD and the constant coefficients as the linear equations in the field $\mathcal{F}_{-1}$ (see Section \ref{sec5} for details). In general, the solutions to these equations are elements of  $\mathcal{F}_{-1}$, i.e., the hyper-functions. However,  some  elements of $\mathcal{F}_{-1}$ can be interpreted as the elements of the ring  $\mathcal{R}_{-1}$, i.e., as the conventional functions (see the embedding \eqref{emb}). In particular, this is valid for all proper rational functions in $S_\kappa$. For the first time, the operational relations of this kind were derived in \cite{Luc21c}. In the rest of this section,  we discuss some important classes of the elements of the field $\mathcal{F}_{-1}$ that can be interpreted as the elements of the ring  $\mathcal{R}_{-1}$ in terms of the so-called convolution series.

The convolution series are a far reaching generalization of the power series. They involve the convolution powers instead of the power law functions:
\begin{equation}
\label{conser}
\Sigma_\kappa(t) = \sum^{+\infty }_{n=0}b_{n}\, \kappa^{<n>}(t),\ b_n \in \Com,\ t>0,
\end{equation}
where $\kappa$ is a Sonin kernel from the class 
$\mathcal{L}_{1}$ (i.e., it belongs to the space $C_{-1}(0,+\infty)$) and the power  series
\begin{equation}
\label{ser}
\Sigma(z) = \sum^{+\infty }_{n=0}b_{n}\, z^n,\ b_{n}\in \Com,\ z\in \Com
\end{equation}
has a non-zero convergence radius. 

As shown in \cite{Luc21c},  the convolution series \eqref{conser} is convergent for any $t>0$. Moreover, its sum $\Sigma_\kappa$ belongs to the space $C_{-1}(0,+\infty)$ and can be interpreted as an  element  of the ring $\mathcal{R}_{-1}$.

Thus, any convergent power series and any Sonin kernel from the class $\mathcal{L}_{1}$ (in fact, any function from the space  $C_{-1}(0,+\infty)$) generate  a convergent convolution series. 

One of the most important for our aims convolution series
\begin{equation}
\label{l}
l_{\kappa,\lambda}(t) := \sum_{n=1}^{+\infty} \lambda^{n-1}\kappa^{<n>}(t),\ \lambda \in \Com,\ t>0
\end{equation}
is generated by the geometric series
\begin{equation}
\label{geom}
\Sigma(z) = \sum_{n=1}^{+\infty} \lambda^{n-1}z^n,\ \lambda \in \Com,\ z\in \Com.
\end{equation}

The series \eqref{l} is convergent for all $t>0$ and the inclusion $l_{\kappa,\lambda} \in C_{-1}(0,+\infty)$ holds true because the convergence radius of the  geometric series is non-zero. 

Let us now consider some examples of the convolution series \eqref{l} generated by different Sonin kernels. We start with kernel $\kappa(t)\equiv 1$ that is  used in the framework of the  Mikusi\'nski operational calculus for the 1st order derivative.  The formula \eqref{alpha_beta} easily implicates the relation $\{ 1\}^{<n>}(t) =  h_1^{<n>}(t) =h_{n}(t)$. The convolution series \eqref{l} takes then the form
\begin{equation}
\label{l-Mic}
l_{\kappa,\lambda}(t) = \sum_{n=1}^{+\infty} \lambda^{n-1}h_n(t) =
\sum_{n=0}^{+\infty} \frac{(\lambda\, t)^n}{n!} = e^{\lambda\, t}.
\end{equation}

In the framework of the operational calculi for the Riemann-Liouville  and for  the Caputo fractional derivatives (see   \cite{LucSri95} and \cite{LucGor99}, respectively), the Sonin kernel $\kappa(t) = h_{\alpha}(t)$ is heavily employed. Using the  formula \eqref{alpha_beta}, we  get the relation $ h_{\alpha}^{<n>}(t) = h_{\alpha\, n}(t)$ and the convolution series \eqref{l} takes the form
\begin{equation}
\label{l-Cap}
l_{\kappa,\lambda}(t) = \sum_{n=1}^{+\infty} \lambda^{n-1}h_{\alpha\, n}(t) =
t^{\alpha-1}\sum_{n=0}^{+\infty} \frac{\lambda^n\, t^{\alpha\, n}}{\Gamma(\alpha\, n+\alpha)} = t^{\alpha -1}E_{\alpha,\alpha}(\lambda\, t^{\alpha}),
\end{equation}
where the two-parameters Mittag-Leffler function $E_{\alpha,\beta}$ is defined by the convergent series    
\begin{equation}
\label{ML}
E_{\alpha,\beta}(z) := \sum_{n=0}^{+\infty} \frac{z^n}{\Gamma(\alpha\, n + \beta) },\ z\in \Com ,\ \beta\in \Com,\ \alpha>0.
\end{equation} 

For other particular cases of the convolution series of type \eqref{l} see \cite{Luc21c}.  

The role of the convolution series $l_{\kappa,\lambda}$ defined by \eqref{l} in operational calculus for the 1st level GFD is illustrated by the following important operational relation derived in \cite{Luc21c} for the first time: 
\begin{equation}
\label{op-rel}
(S_{\kappa} - \lambda)^{-1}\, = \, l_{\kappa,\lambda}(t),\ \lambda \in \Com,\ t>0,
\end{equation}
where by $e^{-1}$ we denote the inverse element to  $e \in \mathcal{F}_{-1}$  with respect to multiplication and the convolution series $l_{\kappa,\lambda} \in \mathcal{R}_{-1}$ is interpreted as an element of $\mathcal{F}_{-1}$.  

The operational relation \eqref{op-rel} holds valid for any Sonin kernel $\kappa \in \mathcal{L}_1$. In particular, the formulas \eqref{l-Mic} and \eqref{l-Cap} for the kernels $\kappa(t) \equiv  1$ and $\kappa(t) = h_{\alpha}(t)$ lead to the following known operational relations: 
\begin{equation}
\label{l-Mic-op}
({S_{\kappa} - \lambda})^{-1}\, = \,   e^{\lambda\, t},\ \kappa(t) \equiv  1,
\end{equation}
\begin{equation}
\label{l-Cap-op}
({S_{\kappa} - \lambda})^{-1}\, = \, t^{\alpha -1}E_{\alpha,\alpha}(\lambda\, t^{\alpha}),\ \kappa(t) = h_{\alpha}(t).
\end{equation}

The operational relation \eqref{op-rel} implicates another important operational relation in terms of the convolution powers of the function  $l_{\kappa,\lambda}\in \mathcal{R}_{-1}$ (\cite{Luc21c}): 
\begin{equation}
\label{op-rel-m}
\left((S_{\kappa} - \lambda)^{-1}\right)^m\, = \, l^{<m>}_{\kappa,\lambda}(t),\ t>0,\ m\in \N.
\end{equation}
The convolution powers $l^{<m>}_{\kappa,\lambda}$ can be calculated in explicit form by constructing the Cauchy products for the convolution series $l_{\kappa,\lambda}$, see \cite{Luc21c} for details and particular cases. 

For the important Sonin kernels $\kappa(t)\equiv 1$ and $\kappa(t) = h_\alpha (t)$, the operational relation \eqref{op-rel-m} takes the following form, respectively:
\begin{equation}
\label{l-Mic-op-m}
\left((S_{\kappa} - \lambda)^{-1}\right)^m\, = \,  h_m(t)\, e^{\lambda\, t},\ \kappa(t)\equiv 1,
\end{equation}
\begin{equation}
\label{l-Cap-op-m}
\left((S_{\kappa} - \lambda)^{-1}\right)^m\, = \, t^{\alpha\, m -1}E_{\alpha,\alpha\, m}^m(\lambda t^\alpha),\ t>0,\ m\in \N,\ \kappa(t) = h_\alpha (t),
\end{equation}
where the Mittag-Leffler type function $E_{\alpha,\beta}^{m}$ is defined by the convergent series
$$
E_{\alpha,\beta}^{m}(z):=
\sum^{\infty
}_{n=0}\frac{(m)_n z^{n}}{  n!
\Gamma (\alpha n + \beta)}, \ \alpha >0, \ \beta \in \Com , \ z\in \Com ,  \
(m)_n = \prod_{i=0}^{n-1} (m+i).
$$

Finally, we mention that the operational relation \eqref{op-rel-m} allows an interpretation of any proper rational function in $S_\kappa$ as an element from the ring $\mathcal{R}_{-1}$. Indeed, let
$$
R(z) = \frac{P(z)}{Q(z)},\ \mbox{deg}(P) < \mbox{deg}(Q),\ z \in \Com 
$$
be a proper rational function with the real or complex coefficients and 
$$
R(z) = \sum_{i=1}^n \sum_{j=1}^{m_i} \frac{a_{ij}}{(z - \lambda_i)^j}
$$
be its 	partial fraction decomposition. 

In the field $\mathcal{F}_{-1}$,  the proper  rational function $R(S_\kappa)$ can be then represented as follows:
\begin{equation}
\label{rat1}
R(S_\kappa) = \frac{P(S_\kappa)}{Q(S_\kappa)} =
\sum_{i=1}^n \sum_{j=1}^{m_i} a_{ij}{\left((S_\kappa - \lambda_i)^{-1}\right)^j}.
\end{equation}

The last formula along with the operational relation \eqref{op-rel-m} leads to an interpretation of the  proper  rational function $R(S_\kappa)$ as an element of the ring $\mathcal{R}_{-1}$ in terms of the convolution series  $l_{\kappa,\lambda}$ and its convolution powers:
\begin{equation}
\label{rat2}
R(S_\kappa) = 
\sum_{i=1}^n \sum_{j=1}^{m_i} a_{ij} l^{<j>}_{\kappa,\lambda_i}(t).
\end{equation}

The operational relation \eqref{rat2} will be used in the next section while applying the operational method for solving the multi-term fractional differential equations with the $n$-fold sequential 1st level GFDs. 



\section{Fractional differential equations with the 1st level GFDs}
\label{sec5}

In this section, we apply the operational calculus developed in the previous section for derivation of the closed form formulas for  solutions to  the linear fractional differential equations with the 1st level GFDs and the $n$-fold sequential 1st level GFDs and with the constant coefficients. The fractional differential equations of this type are new objects not yet considered in the literature. 

To illustrate the method, we start with the simplest equation of this kind
\begin{equation}
\label{eq-0-1}
(  _{1L}\D_{(k_1,k_2)}\, y)(t) = f(t), \ t>0,
\end{equation}
where $ _{1L}\D_{(k_1,k_2)}$ is the 1st level GFD, the given function $f$ belongs to the space $C_{-1}(0,+\infty)$, and the unknown function $y$ is looked for in the space  $ C_{-1,(k_2)}^{1}(0,+\infty)$ defined as in \eqref{C1(k)}. 

According to Theorem \ref{t9}, in the filed $\mathcal{F}_{-1}$ of convolution quotients, the fractional differential equation \eqref{eq-0-1} takes the form of a linear equation 
\begin{equation}
\label{eq-0-2}
S_\kappa \cdot y -  (\I_{(k_2)}\, y)(0)\, k_1  = f,
\end{equation}
where the functions $y,\, k_1,\,  f \in \mathcal{R}_{-1}$ are interpreted as elements of $\mathcal{F}_{-1}$.

Because we did not fix any initial conditions on the unknown function $y$, the term $(\I_{(k_2)}\, y)(0)$ in the equation \eqref{eq-0-2} can be considered to be an arbitrary constant $C\in \R$.

The unique solution to the linear equation \eqref{eq-0-2} has the form
\begin{equation}
\label{eq-0-3}
y = S_{\kappa}^{-1}\cdot f \, +\, S_{\kappa}^{-1} \cdot  (C\, k_1).
\end{equation}

Because of the operational relation $S_{\kappa}^{-1} = \kappa$ (see the formula \eqref{inverse}) and using the embedding of the ring  $\mathcal{R}_{-1}$ into the field  $\mathcal{F}_{-1}$, we can represent the right-hand side of the formula \eqref{eq-0-3} as an element of the ring $\mathcal{R}_{-1}$, i.e., as a conventional function. Thus, the general solution to the equation \eqref{eq-0-1} takes the form
\begin{equation}
\label{eq-0-4}
y(t) = (\kappa \, * f)(t) \, +\, C (\kappa \, * \, k_1)(t), 
\end{equation}
where $C \in \R$ is an arbitrary constant. 

In our derivations, we denoted  the value $(\I_{(k_2)}\, y)(0)$ by the constant $C$. The formula \eqref{eq-0-4} leads then to the statement that the 
 unique solution to the initial-value problem for the fractional differential equation with the 1st level GFD 
\begin{equation}
\label{eq-0-5}
\begin{cases}
(  _{1L}\D_{(k_1,k_2)}\, y)(t) = f(t),\  t>0, \\
(\I_{(k_2)}\, y)(0) = y_0, \ y_0 \in \R 
\end{cases}
\end{equation}
is given by the formula 
\begin{equation}
\label{eq-0-6}
y(t) = (\kappa \, * f)(t) \, +\, y_0 (\kappa \, * \, k_1)(t). 
\end{equation}

The operational method illustrated above can be applied for derivation of the closed form formulas for solutions to other linear fractional differential equations with the 1st level GFDs or to the initial-value problems for such equations. 

\begin{Theorem}
\label{t-de1}
The initial-value problem for the fractional relaxation equation with the 1st level GFD
\begin{equation}
\label{eq-1-1}
\begin{cases}
( _{1L}\D_{(k_1,k_2)}\, y)(t) - \lambda y(t) = f(t), & \lambda \in \R,\ t>0, \\
(\I_{(k_2)}\, y)(0) = y_0, & y_0 \in \R,
\end{cases}
\end{equation}
where the given function $f$ is from the space $C_{-1}(0,+\infty)$ and the unknown function $y$ is looked for in the space  $ C_{-1,(k_2)}^{1}(0,+\infty)$ defined as in \eqref{C1(k)}, possesses a unique solution in the form  
\begin{equation}
\label{eq-1-4}
y(t) = (l_{\kappa,\lambda} \, * f)(t) \, +\, y_0 (l_{\kappa,\lambda} \, * \, k_1)(t), 
\end{equation}
where the convolution series $l_{\kappa,\lambda} \in \mathcal{R}_{-1}$ is defined as in \eqref{l}.
\end{Theorem}

\begin{proof}
To solve the problem \eqref{eq-1-1}, we proceed in the exactly same way as in the case of the equation \eqref{eq-0-1}. First we embed this equation into the field $\mathcal{F}_{-1}$ by using Theorem \ref{t9}:
\begin{equation}
\label{eq-1-2}
S_\kappa \cdot y -  y_0\, k_1 - \lambda y  = f,
\end{equation}
where the functions $y,\, k_1,\, f \in \mathcal{R}_{-1}$ are interpreted as elements of $\mathcal{F}_{-1}$.

The unique solution to the linear equation \eqref{eq-1-2} has the form
\begin{equation}
\label{eq-1-3}
y = (S_{\kappa}-\lambda)^{-1}\cdot f \, +\, y_0\, (S_{\kappa}-\lambda)^{-1} \cdot k_1.
\end{equation}

To represent the hyper-function at the right-hand side of the formula \eqref{eq-1-3} as a conventional function, we employ the operational relation \eqref{op-rel}. Thus we arrive at the solution formula \eqref{eq-1-4} to the initial-value problem \eqref{eq-1-1} and the proof is completed. 
\end{proof} 

The initial-value problem for the fractional relaxation equation \eqref{eq-1-1} with the GFD of the Riemann-Liouville type ($k_1 = h_0$ in the definition of the 1st level GFD) and with the regularized GFD ($k_2 = h_0$ in the definition of the 1st level GFD) were treated for the first time in \cite{Luc22a}  and \cite{Luc21c}, respectively. These are two important particular cases of the problem  \eqref{eq-1-1} worth for being written in explicit form.

\begin{Example}
The initial-value problem  \eqref{eq-1-1} with $k_1 = h_0 \mapsto I,\ k_2 = k$ takes the form
\begin{equation}
\label{eq-1-1-RL}
\begin{cases}
( \D_{(k)}\, y)(t) - \lambda y(t) = f(t), & \lambda \in \R,\ t>0, \\
(\I_{(k)}\, y)(0) = y_0, & y_0 \in \R ,
\end{cases}
\end{equation}
where $\D_{(k)}$ is the  GFD of the Riemann-Liouville type defined by \eqref{1l_RLD}. According to Theorem \ref{t-de1}, it 
has a unique solution $y\in C_{-1,(k)}^{1}(0,+\infty)$ in the form
\begin{equation}
\label{eq-1-4-RL}
y(t) = (l_{\kappa,\lambda} \, * f)(t) \, +\, y_0 l_{\kappa,\lambda}(t). 
\end{equation}

In the case of the Sonin kernels $\kappa(t) = h_\alpha(t),\ k(t) = h_{1-\alpha}(t),\ 0< \alpha <1$, the representation \eqref{l-Cap} leads to the well-known solution formula (see, e.g., \cite{LucSri95})
\begin{equation}
\label{eq-1-4-RL-1}
y(t) = (\tau^{\alpha -1}E_{\alpha,\alpha}(\lambda\, \tau^{\alpha}) \, * f)(t) \, +\, y_0 t^{\alpha -1}E_{\alpha,\alpha}(\lambda\, t^{\alpha}) 
\end{equation}
for the initial-value problem for the fractional differential equation 
\begin{equation}
\label{eq-1-1-RL-1}
\begin{cases}
( D_{0+}^\alpha\, y)(t) - \lambda y(t) = f(t), & \lambda \in \R,\ 0<\alpha <1,\ t>0, \\
(I_{0+}^{1-\alpha}\, y)(0) = y_0, & y_0 \in \R
\end{cases}
\end{equation}
with the Riemann-Liouville fractional derivative $D_{0+}^\alpha$.
\end{Example}

\begin{Example}
Now we consider the initial-value problem  \eqref{eq-1-1} with $k_1 = k,\ k_2 = h_0 \mapsto I$ that can be represented in  the form  
\begin{equation}
\label{eq-1-1-C}
\begin{cases}
( _*\D_{(k)}\, y)(t) - \lambda y(t) = f(t), & \lambda \in \R,\ t>0, \\
y(0) = y_0, & y_0 \in \R ,
\end{cases}
\end{equation}
where ${ } _*\D_{(k)}$ is the regularized GFD defined by \eqref{1l_CD}. 
It's unique solution $y\in C_{-1}^{1}(0,+\infty) :=\{ f\in C_{-1}(0,+\infty): \, f^\prime \in C_{-1}(0,+\infty)\}$ takes the form
\begin{equation}
\label{eq-1-4-C}
y(t) = (l_{\kappa,\lambda} \, * f)(t) \, +\, y_0 L_{\kappa,\lambda}(t), 
\end{equation}
where the function $L_{\kappa,\lambda} \in \mathcal{R}_{-1}$ is defined in form of the following convolution series:
\begin{equation}
\label{L}
L_{\kappa,\lambda}(t) = (l_{\kappa,\lambda}\, *\, k)(t) = 1 + \{ 1 \} * \sum_{n=1}^{+\infty} \lambda^{n}\kappa^{<n>}(t).
\end{equation}

For the Sonin kernels $\kappa(t) = h_{\alpha}(t),\ k(t) = h_{1-\alpha},\ 0<\alpha <1$, the convolution series $l_{\kappa,\lambda}$ is given by the formula \eqref{l-Cap} and the function $L_{\kappa,\lambda}$ can be expressed in terms of the Mittag-Leffler function as follows:
$$
L_{\kappa,\lambda}(t)\, = \, 1 + \{ 1 \} * \sum_{n=1}^{+\infty} \lambda^{n}h_{\alpha\, n}(t) =
1 + \sum_{n=1}^{+\infty} \lambda^{n}h_{\alpha\, n+1}(t) =
$$
\begin{equation}
\label{L-Cap-op}
 \sum_{n=1}^{+\infty} \frac{(\lambda\, t^{\alpha})}{\Gamma (\alpha\, n +1)} = E_{\alpha,1}(\lambda\, t^{\alpha}),
 \end{equation}
where the Mittag-Leffler function $E_{\alpha,1}$ is defined as in  \eqref{ML}.

According to the formula \eqref{eq-1-4-C}, the initial-value problem for the fractional differential equation
\begin{equation}
\label{eq-1-1-2-C}
\begin{cases}
( _*D^\alpha_{0+}\, y)(t) - \lambda y(t) = f(t), & \lambda \in \R,\ 0<\alpha < 1, \ t>0, \\
y(0) = y_0, & y_0\in \R
\end{cases}
\end{equation}
with the Caputo fractional derivative $ _*D^\alpha_{0+}$  has the unique solution in the well-known form (\cite{LucGor99}):
\begin{equation}
\label{sol-1-2}
y (t) = (\tau^{\alpha -1}E_{\alpha,\alpha}(\lambda\, \tau^{\alpha}) \, *\, f)(t) + y_0\, E_{\alpha,1}(\lambda\, t^{\alpha}).
\end{equation}
\end{Example}

Finally, we consider an initial-value problem for the linear multi-term fractional differential equation with the $n$-fold sequential 1st level GFDs and the constant coefficients. The main result is formulated in the next theorem.

\begin{Theorem}
 \label{t-de2}   
The initial-value problem
\begin{equation}
\label{eq-2-1}
\begin{cases}
\sum_{n=0}^{m} b_n(  _{1L}\D_{(k_1,k_2)}^{<n>}\, y)(t) = f(t), \  b_n\in \R,\ b_m \not = 0,\, t>0, \\
(\I_{(k_2)}\,  _{1L}\D_{(k_1,k_2)}^{<n>}\, y)(0) = c_n, \ n=0,\dots,m-1, \ c_n \in \R,
\end{cases}
\end{equation}
 where $ _{1L}\D_{(k_1,k_2)}^{<n>}$ is the $n$-fold sequential 1st level GFD with the kernels from the class $ \mathcal{L}_{1}^{1}$, the given function $f$ is from the space $C_{-1}(0,+\infty)$, and the unknown function $y$ is looked for in the space $C_{-1,(k_1,k_2)}^{m-1}(0,+\infty)$ defined as in \eqref{space}, possesses a unique solution in the form
 \begin{equation}
\label{y-f-4}
y(t) = (f\, * \, G_{\kappa})(t) + (k_1 \, * \, U)(t),
\end{equation}
where  the functions $G_\kappa$ and $U$ are provided in terms of the convolution series $l_{\kappa,\lambda}$ and its convolution powers 
\begin{equation}
\label{G_U}
G(t) = 
\sum_{i=1}^{p} \sum_{j=1}^{p_i} a_{ij} l^{<j>}_{\kappa,\lambda_i}(t),\  
U(t) = \sum_{i=1}^{q} \sum_{j=1}^{q_i} d_{ij} l^{<j>}_{\kappa,\sigma_i}(t).
\end{equation}
The coefficients in the formula \eqref{G_U} are determined by  the partial fractions decompositions of the proper rational functions $\left( P_m(S_\kappa)\right)^{-1}$ and $\frac{Q_{m-1}(S_\kappa)}{P_m(S_\kappa)}$:
\begin{equation}
\label{rat2_1}
\left( P_m(S_\kappa)\right)^{-1} = 
\sum_{i=1}^{p} \sum_{j=1}^{p_i} a_{ij}{\left((S_\kappa - \lambda_i)^{-1}\right)^j},
\end{equation}
\begin{equation}
\label{rat2_2} 
\frac{Q_{m-1}(S_\kappa)}{P_m(S_\kappa)} = \sum_{i=1}^{q} \sum_{j=1}^{q_i} d_{ij}{\left((S_\kappa - \sigma_i)^{-1}\right)^j}.
\end{equation}
 \end{Theorem}
 
\begin{proof}
To solve the initial-value problem \eqref{eq-2-1}, we proceed as in the case of the previous problem \eqref{eq-1-1} with the only difference that we employ Theorem \ref{tsft-n-a} instead of Theorem \ref{t9} and the operational relation \eqref{op-rel-m} instead of the operational relation \eqref{op-rel}. 

By using Theorem \ref{tsft-n-a}, the initial-value problem \eqref{eq-2-1} can be rewritten as a linear equation in the  field $\mathcal{F}_{-1}$ of convolution quotients:
\begin{equation}
\label{eq-2-2}
b_0\, y + \sum_{n=1}^{m} b_n \left( S_\kappa^n \cdot y -  
 \sum_{j=0}^{n-1} c_j \, S_\kappa^{n-j-1} \cdot k_1\right)  = f.
\end{equation}

For the sake of convenience, we introduce the notations
\begin{equation}
\label{poly}
P_m(S_\kappa) = \sum_{n=0}^{m} b_n S_\kappa^n,\ Q_{m-1}(S_\kappa) = \sum_{n=1}^{m} b_n \left(
 \sum_{j=0}^{n-1} c_j \, S_\kappa^{n-j-1}\right).
\end{equation}

With these notations, the linear equation \eqref{eq-2-2} can be rewritten in the form
\begin{equation}
\label{eq-2-3}
P_m(S_\kappa) \cdot y = f +  Q_{m-1}(S_\kappa) \cdot k_1. 
\end{equation}

Its unique solution in the field $\mathcal{F}_{-1}$ (generalized solution to the initial-value problem \eqref{eq-2-3}) is given by the formula
\begin{equation}
\label{eq-2-4}
y = \left( P_m(S_\kappa)\right)^{-1} \cdot f +\frac{Q_{m-1}(S_\kappa)}{P_m(S_\kappa)} \cdot k_1. 
\end{equation}

Finally, we employ the decompositions \eqref{rat2_1} and \eqref{rat2_2}  of the proper rational functions $\left( P_m(S_\kappa)\right)^{-1}$ and $\frac{Q_{m-1}(S_\kappa)}{P_m(S_\kappa)}$, respectively,  and the operational relation \eqref{op-rel-m} to represent the generalized solution \eqref{eq-2-4} as a conventional function from the ring $\mathcal{R}_{-1}$ in form \eqref{y-f-4} and the proof is completed. 
\end{proof}

It is worth mentioning that the part $(f\, * \, G_{\kappa})(t)$ of the solution \eqref{y-f-4} corresponds to the case of the inhomogeneous fractional differential equation and homogeneous initial conditions in the initial-value problem \eqref{eq-2-1} whereas the part $(k_1 \, * \, U)(t)$ is the unique solution to the problem \eqref{eq-2-1} with the homogeneous fractional differential equation and inhomogeneous initial conditions. 

The initial-value problem \eqref{eq-2-1}  with the GFD of the Riemann-Liouville type ($k_1 = h_0$ in the definition of the 1st level GFD) and with the regularized GFD ($k_2 = h_0$ in the definition of the 1st level GFD) was treated for the first time in \cite{Luc22a}  and \cite{Luc21c}, respectively.

\section*{Acknowledgments}
The authors acknowledge the support of the Kuwait University for their joint research project No. SM05/23 "Operational Calculus for the 1st level general fractional derivatives".

\end{document}